\documentclass{svmult}

\usepackage{amsmath}
\usepackage{amsfonts}
\usepackage{amssymb}
\usepackage{hyperref}

\renewcommand{\Re}{\mathrm{Re}}%

\providecommand{\Ent}[1]{\lfloor #1 \rfloor}

\begin{document}

\title*{Asymptotic Approximations to Truncation Errors of Series
  Representations for Special Functions}

\titlerunning{Asymptotic Approximations to Truncation Errors}

\author{Ernst Joachim Weniger}

\institute{Institut f\"{u}r Physikalische und Theoretische Chemie \\
  Universit\"{a}t Regensburg, D-93040 Regensburg, Germany \\
\texttt{joachim.weniger@chemie.uni-regensburg.de}}

\maketitle

\begin{abstract}
\noindent
Asymptotic approximations ($n \to \infty$) to the truncation errors $r_n
= - \sum_{\nu=0}^{\infty} a_{\nu}$ of infinite series
$\sum_{\nu=0}^{\infty} a_{\nu}$ for special functions are constructed by
solving a system of linear equations. The linear equations follow from an
approximative solution of the inhomogeneous difference equation $\Delta
r_n = a_{n+1}$. In the case of the remainder of the Dirichlet series for
the Riemann zeta function, the linear equations can be solved in closed
form, reproducing the corresponding Euler-Maclaurin formula.  In the case
of the other series considered -- the Gaussian hypergeometric series
${}_2 F_1 (a, b; c; z)$ and the divergent asymptotic inverse power series
for the exponential integral $E_1 (z)$ -- the corresponding linear
equations are solved symbolically with the help of Maple. The
practical usefulness of the new formalism is demonstrated by some
numerical examples.
\end{abstract}

\typeout{==> Section 1: Introduction}
\section{Introduction}
\label{Sec:Intro}

A large part of special function theory had been developed already in the
19th century. Thus, it is tempting to believe that our knowledge about
special functions is essentially complete and that no significant new
developments are to be expected. However, up to the middle of the 20th
century, research on special functions had emphasized analytical results,
whereas the efficient and reliable evaluation of most special functions
had been -- and to some extend still is -- a more or less unsolved
problem.

Due to the impact of computers on mathematics, the situation has changed
substantially. We witness a revival of interest in special functions.
The general availability of electronic computers in combination with the
development of powerful computer algebra systems like Maple or
Mathematica opened up many new applications, and it also created a great
demand for efficient and reliable computational schemes (see for example
\cite{Gil/Segura/Temme/2003,Lozier/Olver/1994,VanDerLaan/Temme/1980} or
\cite[Section 13]{Temme/1996a} and references therein).

Most special functions are defined via infinite series. Examples are the
Dirichlet series for the Riemann zeta function,
\begin{equation}
  \label{DiriSerZetaFun}
\zeta (s) \; = \; \sum _{\nu=0}^{\infty} \; (\nu+1)^{-s} \, .
\end{equation}
which converges for $\Re (s) > 1$, or the Gaussian hypergeometric series
\begin{equation}
  \label{Ser_2F1}
{}_2 F_1 (a, b; c; z) \; = \; \sum_{\nu=0}^{\infty} \,
\frac {(a)_{\nu} (b)_{\nu}} {(c)_{\nu} {\nu}!} \, z^{\nu} \, ,
\end{equation}
which converges for $\vert z \vert < 1$.  

The definition of special functions via infinite series is to some extent
highly advantageous since it greatly facilitates analytical
manipulations. However, from a purely numerical point of view, infinite
series representations are at best a mixed blessing. For example, the
Dirichlet series (\ref{DiriSerZetaFun}) converges for $\Re (s) > 1$, but
is notorious for extremely slow convergence if $\Re (s)$ is only slightly
larger than one. Similarly, the Gaussian hypergeometric series
(\ref{Ser_2F1}) converges only for $\vert z \vert < 1$, but the
corresponding Gaussian hypergeometric function is a multivalued function
defined in the whole complex plane with branch points at $z = 1$ and
$\infty$. A different computational problem occurs in the case of the
asymptotic series for the exponential integral:
\begin{equation}
\label{AsySerE_1}
z \, \mathrm{e}^z \, E_1 (z) \, \sim \,
\sum_{m=0}^{\infty} \, (-1/z)^m \, m! \; = \;
{}_2 F_0 (1, 1; - 1/z) \, , \qquad z \to \infty \, .
\end{equation}
This series is probably the most simple example of a large class of
series that diverge for every finite argument $z$ and that are only
asymptotic in the sense of Poincar\'{e} as $z \to \infty$. In contrast,
the exponential integral $E_1 (z)$, which has a cut along the negative
real axis, is defined in the whole complex plane.

Problems with slow convergence or divergence were encountered already in
the early days of calculus. Thus, numerical techniques for the
acceleration of convergence or the summation of divergent series are
almost as old as calculus. According to Knopp \cite[p.\ 
249]{Knopp/1964}, the first systematic work in this direction can be
found in Stirling's book \cite{Stirling/1730}, which was published
already in 1730 (recently, Tweddle \cite{Tweddle/2003} published a new
annotated translation), and in 1755 Euler \cite{Euler/1755} published
the series transformation which now bears his name.

The convergence and divergence problems mentioned above can be
formalized as follows: Let us assume that the partial sums $s_n =
\sum_{\nu=0}^{n} a_{\nu}$ of a convergent or divergent but summable
series form a sequence $\{ s_n \}_{n=0}^{\infty}$ whose elements can be
partitioned into a (generalized) limit $s$ and a remainder or truncation
error $r_n$ according to
\begin{equation}
  \label{s_n_r_n}
s_n \; = \; s \, + \, r_n \, , \qquad n \in \mathbb{N}_0 \, .
\end{equation}
This implies
\begin{equation}
  \label{DefRemInfSer}
r_n \; = \; - \sum_{\nu=n+1}^{\infty}\, a_{\nu} \, , 
\qquad n \in \mathbb{N}_0 \, .
\end{equation}

At least in principle, a convergent infinite series can be evaluated by
adding up the terms successively until the remainders become negligible.
This approach has two obvious shortcomings. Firstly, convergence can be
so slow that it is uneconomical or practically impossible to achieve
sufficient accuracy. Secondly, this approach does not work in the case of
a divergent but summable series because increasing the index $n$ normally
only aggravates divergence.
 
As a principal alternative, we can try to compute a sufficiently accurate
approximation $\bar{r}_n$ to the truncation error $r_n$. If this is
possible, $\bar{r}_n$ can be eliminated from $s_n$, yielding a (much)
better approximation $s_n-\bar{r}_n$ to the (generalized) limit $s$ than
$s_n$ itself.

This approach looks very appealing since it is in principle remarkably
powerful. In addition, it can avoid the troublesome asymptotic regime of
large indices $n$, and it also works in the case of divergent but
summable sequences and series. Unfortunately, it is by no means easy to
obtain sufficiently accurate approximations $\bar{r}_n$ to truncation
errors $r_n$. The straightforward computation of $r_n$ by adding up the
terms does not gain anything.

The Euler-Maclaurin formula, which is discussed in Section
\ref{Sec:EulerMac}, is a principal analytical tool that produces
asymptotic approximations to truncation errors of monotone series in
terms of integrals plus correction terms. Unfortunately, it is not always
possible to apply the Euler-Maclaurin formula. Given a reasonably well
behaved integrand, it is straightforward to compute a sum of integrand
values plus derivatives of the integrand. But for a given series term
$a_n$, it may be prohibitively difficult to differentiate and integrate
it with respect to the index $n$.

In Section \ref{Sec:AsyTruncErrAppr}, an alternative approach for the
construction of asymptotic approximations ($n \to \infty$) to the
truncation errors $r_n$ of infinite series is proposed that is based on
the solution of a system of linear equations. The linear equations exist
under very mild conditions: It is only necessary that the ratio
$a_{n+2}/a_{n+1}$ or similar ratios of series terms possesses an
asymptotic expansion in terms of inverse powers $1/(n+\alpha)$ with
$\alpha > 0$. Moreover, it is also fairly to solve these linear
equations since they have a triangular structure.

The asymptotic nature of the approximants makes it difficult to use them
also for small indices $n$, although this would be highly desirable. In
Section \ref{Sec:NumAnaCon}, it is mentioned briefly that factorial
series and Pad\'{e} approximants can be helpful in this respect since
they can accomplish a numerical analytic continuation.

In Section \ref{Sec:DirSerRiemannZetaFun}, the formalism proposed in this
article is applied to the truncation error of the Dirichlet series for
the Riemann zeta function. It is shown that the linear equations can in
this case be reduced to a well known recurrence formula of the Bernoulli
numbers. Accordingly, the terms of the corresponding Euler-Maclaurin
formula are exactly reproduced.

In Section \ref{Sec:Gauss_Hyg_2F1}, the Gaussian hypergeometric series
${}_2 F_1 (a, b; c; z)$ is treated. Since the terms of this series depend
on three parameters and one argument, a closed form solution of the
linear equations seems to be out of reach. Instead, approximations are
computed symbolically with the help of the computer algebra system Maple.
The practical usefulness of these approximations is demonstrated by some
numerical examples.

In Section \ref{Sec:AsySer_E_1}, the divergent asymptotic inverse power
series for the exponential integral $E_1 (z)$ is treated. Again, the
linear equations are solved symbolically with the help of Maple, and the
practical usefulness of these solutions is demonstrated by some numerical
examples.

\typeout{==> Section 2: The Euler-Maclaurin Formula}
\section{The Euler-Maclaurin Formula}
\label{Sec:EulerMac}

The derivation of the Euler-Maclaurin formula is based on the assumption
that $g (x)$ is a smooth and slowly varying function. Then, $\int_{M}^{N}
g (x) \mathrm{d} x$ with $M, N \in \mathbb{Z}$ can be approximated by the
finite sum $\frac{1}{2} g (M) + g (M+1) + \dots + g (N-1) + \frac{1}{2} g
(N)$. This finite sum can also be interpreted as a trapezoidal quadrature
rule. In the years between 1730 and 1740, Euler and Maclaurin derived
independently correction terms to this quadrature rule, which ultimately
yielded what we now call the Euler-Maclaurin formula (see for example
\cite[Eq.\ (1.20)]{Temme/1996a}):
\begin{subequations}
  \label{EuMaclau}
  \begin{align}
%  \label{}
    \sum_{\nu=M}^{N} \, g (\nu) & \; = \; \int_{M}^{N} \, g (x) \,
    \mathrm{d} x \, + \, \frac{1}{2} \,
    \bigl[g (M) + g (N)\bigr] \notag \\
    & \qquad + \, \sum_{j=1}^{k} \, \frac {B_{2j}} {(2j)!}  \left[
      g^{(2j-1)} (N) - g^{(2j-1)} (M) \right] \, + \, R_k (g) \, ,
    \\
%  \label{}
    R_k (g) & \; = \; - \, \frac{1}{(2k)!} \, \int_{M}^{N} \, B_{2k}
    \bigl( x - \Ent{x} \bigr) \, g^{(2k)} (x) \, \mathrm{d} x \, .
  \end{align}
\end{subequations}
Here, $g^{(m)} (x)$ is the $m$-th derivative, $\Ent{x}$ is the integral
part of $x$, $B_m (x)$ is a Bernoulli polynomial defined by the
generating function $t \mathrm{e}^{xt}/(\mathrm{e}^t - 1) =
\sum_{n=0}^{\infty} B_n (x) t^n / n!$, and $B_m = B_m (0)$ is a Bernoulli
number. 

It is not a priori clear whether the integral $R_k (g)$ in
(\ref{EuMaclau}) vanishes as $k \to \infty$ for a given function $g (x)$.
Thus, the Euler-Maclaurin formula may lead to an asymptotic expansion
that ultimately diverges. In this article, it is always assumed that the
Euler-Maclaurin formula and related expansions are only asymptotic in the
sense of Poincar\'{e}.

Although originally used to express the in the early 18th century still
unfamiliar integral in terms more elementary quantities, the
Euler-Maclaurin formula is now often used to approximate the truncation
error $r_n = - \sum_{\nu=n+1}^{\infty} a_{\nu}$ of a slowly convergent
monotone series by an integral plus correction terms. The power and the
usefulness of this approach can be demonstrated convincingly via the
Dirichlet series (\ref{DiriSerZetaFun}) for the Riemann zeta function.

The terms $(\nu+1)^{-s}$ of the Dirichlet series (\ref{DiriSerZetaFun})
are obviously smooth and slowly varying functions of the index $\nu$, and
they can be differentiated and integrated easily. Thus, the application
of the Euler-Maclaurin formula (\ref{EuMaclau}) with $M = n+1$ and $N =
\infty$ to the truncation error of the Dirichlet series yields:
\begin{subequations}
  \label{EuMacZeta}
  \begin{align}
  \label{EuMacZeta_a}
  - \, \sum_{\nu=n+1}^{\infty} \, (\nu+1)^{-s} & \; = \;
  - \, \frac{(n+2)^{1-s}}{s-1} \, - \, \frac{1}{2} \, (n+2)^{-s}
  \notag \\
  & \qquad \, - \, \sum_{j=1}^{k} \, \frac {(s)_{2j-1} \, B_{2j}} {(2j)!}
  \,
  (n+2)^{-s-2j+1} \, + \, R_k (n,s) \, , \\
    \label{EuMacZeta_b}
    R_k (n,s) & \; = \; \frac{(s)_{2k}}{(2k)!} \,
    \int_{n+1}^{\infty} \, \frac {B_{2k} \bigl( x - \Ent{x} \bigr)}
    {(x+1)^{s+2k}} \, \mathrm{d} x \, .
   \end{align}
\end{subequations}
Here, $(s)_m = s (s+1) \cdots (s+m-1) = \Gamma(s+m)/\Gamma(s)$ with $s
\in \mathbb{C}$ and $m \in \mathbb{N}_0$ is a Pochhammer symbol.

In \cite[Tables 8.7 and 8.8, p.\ 380]{Bender/Orszag/1978} and in
\cite[Section 2]{Weniger/Kirtman/2003} it was shown that a few terms of
the sum in (\ref{EuMacZeta_a}) suffice for a convenient and reliable
computation of $\zeta (s)$ with $s = 1.1$ and $s = 1.01$, respectively.
For these arguments, the Dirichlet series for $\zeta (s)$ converges so
slowly that it is practically impossible to evaluate it by adding up its
terms.

In order to understand better its nature, the Euler-Maclaurin formula
(\ref{EuMaclau}) is rewritten in a more suggestive form. Let us set $M =
n+1$ and $N = \infty$, and let us also assume $\lim_{N \to \infty} g (N)
= \lim_{N \to \infty} g' (N) = \lim_{N \to \infty} g'' (N) = \dots = 0$.
With the help of $B_0=1$. $B_1=-1/2$, and $B_{2n+1}=0$ with $n \in
\mathbb{N}$ (see for example \cite[p.\ 3]{Temme/1996a}), we obtain:
\begin{subequations}
  \label{EuMaclau_mod1}
  \begin{align}
    \label{EuMaclau_mod1_a}
    - \, \sum_{\nu=n+1}^{\infty} \, g (\nu) & \; = \; 
    - B_0 \, \int_{n+1}^{\infty} \, g (x) \, \mathrm{d} x
    \, + \, \sum_{\mu=1}^{m} \, \frac{(-1)^{\mu-1} B_{\mu}}{\mu!} \, 
    g^{(\mu-1)} (\nu) \, + \, R_m (g) \, ,
    \\
    \label{EuMaclau_mod1_b}
    R_m (g) & \; = \; \frac{(-1)^m}{(m)!} \,
    \int_{n+1}^{\infty} \, B_{m} \bigl( x - \Ent{x} \bigr)
    g^{(m)} (x) \, \mathrm{d} x \, , \quad m \in \mathbb{N} \, .
  \end{align}
\end{subequations}
In the same way, we obtain for the Euler-Maclaurin approximation
(\ref{EuMacZeta}) to the truncation error of the Dirichlet series:
\begin{subequations}
  \label{EuMacZeta_mod2}
  \begin{align}
    \label{EuMacZeta_mod2_a}
    - \, \sum_{\nu=n+1}^{\infty} \, (\nu+1)^{-s} & \; = \;
    \sum_{\mu=0}^{m} \, \frac{(-1)^{\mu-1} (s)_{\mu-1} B_{\mu}}{\mu!} \,
    (n+2)^{1-s-\mu} \, + \, R_m (n,s) \, ,
    \\
    \label{EuMacZeta_mod2_b}
    R_m (n,s) & \; = \; \frac{(-1)^m (s)_{m}}{(m)!} \,
    \int_{n+1}^{\infty} \, \frac {B_{m} \bigl( x - \Ent{x} \bigr)}
    {(1+x)^{s+m}} \, \mathrm{d} x \, , \quad m \in \mathbb{N} \, .
  \end{align}
\end{subequations}

The reformulated Euler-Maclaurin approximation (\ref{EuMacZeta_mod2})
looks suspiciously like a truncated expansion of the truncation error in
terms of the asymptotic sequence $\{ (n+2)^{-\mu} \}_{\mu=0}^{\infty}$
of inverse powers. An analogous interpretation of the reformulated
Euler-Maclaurin formula (\ref{EuMaclau_mod1}) is possible if we assume
that the quantities $\int_{n+1}^{\infty} g (x) \mathrm{d}x, g (n), g'
(n), g'' (n), \dots$ form an asymptotic sequence $\bigl\{
\mathcal{G}_{\mu} (n) \bigr\}_{\mu=0}^{\infty}$ as $n \to \infty$
according to
\begin{subequations}
  \label{AsySeq_EuMaclau_mod}
  \begin{align}
%  \label{AsySeq_EuMaclau_mod_a}
    \mathcal{G}_0 (n) & \; = \; \int_{n+1}^{\infty} \, g (x) \,
    \mathrm{d} x \, ,
    \\
%  \label{AsySeq_EuMaclau_mod_b}
    \mathcal{G}_{\mu} (n) & \; = \; g^{(\mu-1)} (n) \, , \qquad \mu \in
    \mathbb{N} \, .
\end{align}
\end{subequations}

The expansion of the truncation error $-\sum_{\nu=n+1}^{\infty} g (\nu)$
in terms of the asymptotic sequence $\bigl\{ \mathcal{G}_{\mu} (n)
\bigr\}_{\mu=0}^{\infty}$ according to (\ref{EuMaclau_mod1}) has the
undeniable advantage that the expansion coefficients do not depend on
the terms $g (\nu)$ and are explicitly known. The only remaining
computational problem is the determination of the leading elements of
the asymptotic sequence $\bigl\{ \mathcal{G}_{\mu} (n)
\bigr\}_{\mu=0}^{\infty}$. In the case of the Dirichlet series
(\ref{DiriSerZetaFun}), this is trivially simple. Unfortunately,
the terms of most series expansions for special functions are (much)
more complicated than the terms of the Dirichlet series
(\ref{DiriSerZetaFun}). In those less fortunate cases, it can be
extremely difficult to do the necessary differentiations and
integrations. Thus, the construction of the asymptotic sequence $\bigl\{
\mathcal{G}_{\mu} (n) \bigr\}_{\mu=0}^{\infty}$ may turn out to be an
unsurmountable problem.

\typeout{==> Section 3: Asymptotic Approximations To Truncation Errors}
\section{Asymptotic Approximations To Truncation Errors}
\label{Sec:AsyTruncErrAppr}

Let us assume that we want to construct an asymptotic expansion of a
special function $f (z)$ as $z \to \infty$. First, we have to find a
suitable asymptotic sequence $\{ \varphi_{j} (z) \}_{j=0}^{\infty}$.
Obviously, $\{ \varphi_{j} (z) \}_{j=0}^{\infty}$ must be able to model
the essential features of $f (z)$ as $z \to \infty$. On the other hand,
$\{ \varphi_{j} (z) \}_{j=0}^{\infty}$ should also be sufficiently simple
in order to facilitate the necessary analytical manipulations. In that
respect, the most convenient asymptotic sequence is the sequence $\{
z^{-j} \}_{j=0}^{\infty}$ of inverse powers, and it is also the one which
is used almost exclusively in special function theory. An obvious example
is the asymptotic series (\ref{AsySerE_1}).

The behavior of most special functions as $z \to \infty$ is incompatible
with an expansion in terms of inverse powers. Therefore, an indirect
approach has to be pursued: For a given $f (z)$, one has to find some $g
(z)$ such that the ratio $f (z)/g (z)$ admits an asymptotic expansion in
terms of inverse powers:
\begin{equation}
  \label{GenAsySer_SpecFun}
f (z)/g (z) \; \sim \; \sum_{j=0}^{\infty} c_j/z^j \, ,
\qquad z \to \infty \, .
\end{equation}
Although $f (z)$ cannot be expanded in terms of inverse powers $\{ z^{-j}
\}_{j=0}^{\infty}$, it can be expanded in terms of the asymptotic
sequence $\{ g (z)/z^{j} \}_{j=0}^{\infty}$. The asymptotic series
(\ref{AsySerE_1}) is of the form of (\ref{GenAsySer_SpecFun}) with $f (z)
= E_1 (z)$ and $g (z) = \exp (-z)/z$.

It is the central hypothesis of this article that such an indirect
approach is useful for the construction of asymptotic approximations to
remainders $r_n = - \sum_{\nu=n+1}^{\infty} a_{\nu}$ of infinite series
as $n \to \infty$. Thus, instead of trying to use the technically
difficult Euler-Maclaurin formula (\ref{EuMaclau}), we should try to find
some $\rho_n$ such that the ratio $r_n/\rho_n$ admits an asymptotic
expansion as $n \to \infty$ in terms of inverse powers $\{
(n+\alpha)^{-j} \}_{j=0}^{\infty}$ with $\alpha > 0$.

A natural candidate for $\rho_n$ is the first term $a_{n+1}$ neglected in
the partial sum $s_n = \sum_{\nu=0}^{n} a_{\nu}$, but in some cases it is
better to choose instead $\rho_n = a_n$ or $\rho_n = (n+\alpha) a_{n+1}$
with $\alpha > 0$. Moreover, the terms $a_{n+1}$ and the remainders $r_n$
of an infinite series are connected by the inhomogeneous difference
equation
\begin{equation}
  \label{FinDiffEqTerm}
\Delta r_n \; = \; r_{n+1} - r_n \; = \; a_{n+1} \, ,
\qquad n \in \mathbb{N}_0 \, .
\end{equation}
In Jagerman's book \cite[Chapter 3 and 4]{Jagerman/2000}, solutions to
difference equations of that kind are called N\"{o}rlund sums.

If we knew how to solve (\ref{FinDiffEqTerm}) efficiently and reliably
for essentially arbitrary inhomogeneities $a_{n+1}$, all problems
related to the evaluation of infinite series would in principle be
solved. Unfortunately, this is not the case.  Nevertheless, we can use
(\ref{FinDiffEqTerm}) to construct the leading terms of an asymptotic
expansion of $r_n/a_{n+1}$ or of related expressions in terms of inverse
powers.

For that purpose, we make the following ansatz:
\begin{equation}
  \label{Ansatz_RemSer}
r_{n}^{(m)} \; = \; - \, a_{n+1} \, \sum_{\mu=0}^{m} \,
\frac{\gamma_{\mu}^{m}}{(n+\alpha)^{\mu}} \, ,
\qquad n \in \mathbb{N}_0 \, , \quad m \in \mathbb{N} \, , 
\quad \alpha > 0 \, .
\end{equation}
This ansatz, which is inspired by the theory of converging factors
\cite{Airey/1937,Miller/1952} and by a truncation error estimate for
Levin's sequence transformation \cite{Levin/1973} proposed by Smith and
Ford \cite[Eq.\ (2.5)]{Smith/Ford/1979} (see also \cite[Section
7.3]{Weniger/1989} or \cite[Section IV]{Weniger/2004}), is not completely
general and has to be modified slightly both in the case of the Dirichlet
series (\ref{DiriSerZetaFun}) for the Riemann zeta function, which is
discussed in Section \ref{Sec:DirSerRiemannZetaFun}, and in the case of
the divergent asymptotic series (\ref{AsySerE_1}) for the exponential
integral, which is discussed in Section \ref{Sec:AsySer_E_1}. Moreover,
the ansatz (\ref{Ansatz_RemSer}) does not cover the series expansions of
all special functions of interest. For example, in \cite{Weniger/2003} a
power series expansion for the digamma function $\psi (z)$ was analyzed
whose truncation errors cannot be approximated by a truncated power
series of the type of (\ref{Ansatz_RemSer}). Nevertheless, the examples
considered in this article should suffice to convince even a sceptical
reader that the ansatz (\ref{Ansatz_RemSer}) is indeed computationally
useful.

We cannot expect that the ansatz (\ref{Ansatz_RemSer}) satisfies the the
inhomogeneous difference equation (\ref{FinDiffEqTerm}) exactly.
However, we can choose the unspecified coefficients $\gamma_{\mu}^{(m)}$
in (\ref{Ansatz_RemSer}) in such a way that only a higher order error
remains:
\begin{align}
  \label{RatioExpan}
  \frac{r_{n+1}^{(m)} - r_{n}^{(m)}}{a_{n+1}} & \; = \; \sum_{\mu=0}^{m}
  \, \frac{\gamma_{\mu}^{m}}{(n+\alpha)^{\mu}} \, - \,
  \frac{a_{n+2}}{a_{n+1}} \, \sum_{\mu=0}^{m} \,
  \frac{\gamma_{\mu}^{m}}{(n+\alpha+1)^{\mu}}
  \\
  \label{ConvCond}
  & \; = \; 1 + \mathrm{O} (n^{-m-1}) \, , \qquad n \to \infty \, .
\end{align}
The approach of this article depends crucially on the assumption that the
ratio $a_{n+2}/a_{n+1}$ can be expressed as an (asymptotic) power series
in $1/(n+\alpha)$. If this is the case, then the right-hand side of
(\ref{RatioExpan}) can be expanded in powers of $1/(n+\alpha)$ and we
obtain:
\begin{equation}
%  \label{}
\frac{r_{n+1}^{(m)} - r_{n}^{(m)}}{a_{n+1}} \; = \;
\sum_{\mu=0}^{m} \, \frac{\mathcal{C}_{\mu}^{(m)}}{(n+\alpha)^{\mu}} 
\, + \, \mathrm{O} \bigl( n^{-m-1} \bigr) \, , \qquad n \to \infty \, .
\end{equation}
Now, (\ref{ConvCond}) implies that we have solve the following system of
linear equations:
\begin{equation}
  \label{SysLinEqs}
\mathcal{C}_{\mu}^{(m)} \; = \; \delta_{\mu 0} \, , 
\qquad 0 \le \mu \le m \, .
\end{equation}
Since $\mathcal{C}_{\mu}^{(m)}$ with $0 \le \mu \le m$ contains only the
unspecified coefficients $\gamma_{0}^{(m)}$, $\dots$,
$\gamma_{\mu}^{(m)}$ but not $\gamma_{\mu+1}^{(m)}$ $\dots$,
$\gamma_{m}^{(m)}$, the linear system (\ref{SysLinEqs}) has a triangular
structure and the unspecified coefficients $\gamma_{0}^{(m)}, \dots,
\gamma_{\mu}^{(m)}$ can be determined by solving successively the
equations $\mathcal{C}_{0}^{(m)} = 1$, $\mathcal{C}_{1}^{(m)} = 0$,
$\dots$, $\mathcal{C}_{m}^{(m)} = 0$.

Another important aspect is that the linear equations (\ref{SysLinEqs})
do not depend explicitly on $m$, which implies that the coefficients
$\gamma_{\mu}^{(m)}$ in (\ref{Ansatz_RemSer}) also do not depend
explicitly on $m$. Accordingly, the superscript $m$ of both
$\mathcal{C}_{\mu}^{(m)}$ and $\gamma_{\mu}^{(m)}$ is superfluous and
will be dropped in the following Sections.

\typeout{==> Section 4: Numerical Analytic Continuation}
\section{Numerical Analytic Continuation}
\label{Sec:NumAnaCon}

Divergent asymptotic series of the type of (\ref{GenAsySer_SpecFun}) can
be extremely useful computationally: For sufficiently large arguments
$z$, truncated expansions of that kind are able to provide (very)
accurate approximations to the corresponding special functions, in
particular if the series is truncated in the vicinity of the minimal
term. If, however, the argument $z$ is small, truncated expansions of
that kind produce only relatively poor or even completely nonsensical
results.

We can expect that our asymptotic expansions in powers of $1/(n+\alpha)$
have similar properties. Thus, we can be confident that they produce
(very) good results for sufficiently large indices $n$, but it would be
overly optimistic to assume that these expressions necessarily produce
good results in the nonasymptotic regime of moderately large or even
small indices $n$.

Asymptotic approximants can often be constructed (much) more easily than
other approximants that are valid in a wider domain. Thus, it is
desirable to use asymptotic approximants also outside the asymptotic
domain. This means that we would like to use our asymptotic approximants
also for small indices $n$ in order to avoid the computationally
problematic asymptotic regime of large indices.  Obviously, this is
intrinsically contradictory. We also must find a way of extracting
additional information from the terms of a truncated divergent inverse
power series expansion.

Often, this can be accomplished at low computational costa by converting
an inverse power series $\sum_{n=0}^{\infty} c_n/z^n$ to a factorial
series $\sum_{n=0}^{\infty} \tilde{c}_n/(z)_n$. Factorial series, which
had already been known to Stirling \cite[p.\ 6]{Tweddle/2003},
frequently have superior convergence properties. An example is the
incomplete gamma function $\Gamma (a, z)$, which possesses a divergent
asymptotic series of the type of (\ref{GenAsySer_SpecFun}) \cite[Eq.\ 
(6) on p.\ 135]{Erdelyi/Magnus/Oberhettinger/Tricomi/1953b} and also a
convergent factorial series \cite[Eq.\ (1) on p.\ 
139]{Erdelyi/Magnus/Oberhettinger/Tricomi/1953b}. Accordingly, the
otherwise so convenient inverse powers are not necessarily the
computationally most effective asymptotic sequence.

The transformation of an inverse power series to a factorial series can
be accomplished with the help of the Stirling numbers of the first kind
which are normally defined via the expansion $(z-n+1)_n =
\sum_{\nu=0}^{n} \textbf{S}^{(1)} (n, \nu) z^{\nu}$ of a Pochhammer
symbol in terms of powers. As already known to Stirling (see for example
\cite[p.\ 29]{Tweddle/2003} or \cite[Eq.\ (6) on p.\ 78]{Nielsen/1965}),
the Stirling numbers of the first kind occur also in the factorial series
expansion of an inverse power:
\begin{equation}
  \label{GenFunSti1}
\frac{1}{z^{k+1}} \; = \; \sum_{\kappa=0}^{\infty} \,
\frac{(-1)^{\kappa} \, \textbf{S}^{(1)} (k+\kappa, k)}{(z)_{k+\kappa+1}}
\, , \qquad k \in \mathbb{N}_0 \, .
\end{equation}
This infinite generating function can also be derived by exploiting the
well known recurrence relationships of the Stirling numbers.

With the help of (\ref{GenFunSti1}), the following transformation formula
can be derived easily:
\begin{equation}
  \label{PowSerg_2FactSer}
\sum_{n=0}^{\infty} \, \frac{c_n}{z^{n}} \; = \; 
c_0 \, + \, \frac{c_1}{(z)_1}  \, + \,
\sum_{k=2}^{\infty} \, \frac{(-1)^{k}}{(z)_{k}} \,
\sum_{\kappa=1}^{k} \, (-1)^{\kappa} \, \textbf{S}^{(1)} (k-1, \kappa-1)
\, c_{\kappa} \, .
\end{equation}

Let us now assume that the coefficients $\gamma_{\mu}$ with $0 \le \mu
\le m$ of a truncated expansion of $r_n/a_{n+1}$ in powers of
$1/(n+\alpha)$ according to (\ref{Ansatz_RemSer}) are known. Then,
(\ref{PowSerg_2FactSer}) implies that we can use the transformation
scheme
\begin{equation}
%  \label{gamma2Gamma}
\tilde{\gamma_{\mu}} \; = \;
\begin{cases}
  \gamma_{\mu} \, , \qquad \mu = 0, 1 \, , \\[1.5\jot]
  {\displaystyle 
   \sum_{\nu=1}^{\mu} \, (-1)^{\mu+\nu} \, \textbf{S}^{(1)} (\mu-1, \nu-1)
   \, \gamma_{\nu} \, , \quad \mu \ge 2} \, ,
\end{cases}
\end{equation}
to obtain instead of (\ref{Ansatz_RemSer}) the truncated factorial series
\begin{equation}
%  \label{}
\tilde{r}_{n}^{(m)} \; = \; - a_{n+1} \, \left[
\sum_{\mu=0}^{m} \, \frac{\tilde{\gamma_{\mu}}}{(n+\alpha)_{\mu}} 
\, + \, \mathrm{O} \bigl( n^{-m-1} \bigr) \right]  \, ,
\qquad n \to \infty \, .
\end{equation}

Pad\'{e} approximants, which convert the partial sums of a formal power
series to a doubly indexed sequence of rational functions, can also be
quite helpful. They are now used almost routinely in applied mathematics
and theoretical physics to overcome convergence problems with power
series (see for example the monograph by Baker and Graves-Morris
\cite{Baker/Graves-Morris/1996} and references therein).  The ansatz
(\ref{Ansatz_RemSer}) produces a truncated series expansion of
$r_n/a_{n+1}$ in powers of $1/(n+\alpha)$, which can be converted to a
Pad\'{e} approximant, i.e., to a rational function in $1/(n+\alpha)$.

The numerical results presented in Sections \ref{Sec:Gauss_Hyg_2F1} and
\ref{Sec:AsySer_E_1} that the conversion to factorial series and Pad\'{e}
approximants improves the accuracy of our approximants, in particular for
small indices $n$.

\typeout{==> Section 5: The Dirichlet Series for the Riemann Zeta Function}
\section{The Dirichlet Series for the Riemann Zeta Function}
\label{Sec:DirSerRiemannZetaFun}

In this Section, an asymptotic approximation to the truncation error of
the Dirichlet series for the Riemann zeta function is constructed by
suitably adapting the approach described in Section
\ref{Sec:AsyTruncErrAppr}. In the case of the Dirichlet series
(\ref{DiriSerZetaFun}), we have:
\begin{align}
  \label{ParSumDirSerZeta}
  s_n & \; = \; \sum_{\nu=0}^{n} \, (\nu+1)^{-s} \, ,
  \\
  \label{RemDirSerZeta}
  r_n & \; = \; - \sum_{\nu=n+1}^{\infty} \, (\nu+1)^{-s} \; = \; -
  (n+2)^{-s} \, \sum_{\nu=0}^{\infty} \,
  \left(1+\frac{\nu}{n+2}\right)^{-s} \, ,
  \\
  \label{FinDifEq_SerZeta}
  \Delta r_n & \; = \; (n+2)^{-s} \, .
\end{align}

It is an obvious idea to express the infinite series on the right-hand
side of (\ref{RemDirSerZeta}) as a power series in $1/(n+2)$. If $\nu <
n+2$, we can use the binomial series $(1+z)^a = {}_1 F_0 (-a; -z) =
\sum_{m=0}^{\infty} \binom{a}{m} z^m$ \cite[p.\
38]{Magnus/Oberhettinger/Soni/1966}, which converges for $\vert z \vert <
1$. We thus obtain
\begin{equation}
%  \label{}
[1+\nu/(n+2)]^{-s} \; = \; \sum_{m=0}^{\infty} \, \frac{(s)_m}{m} \,
[-\nu/(n+2)]^m \, .
\end{equation}
The infinite series converges if $\nu/(n+2) < 1$. Thus, an expansion of
the right-hand side of (\ref{RemDirSerZeta}) in powers of $1/(n+2)$ can
only be asymptotic as $n \to \infty$. Nevertheless, a suitably truncated
expansion suffices for our purposes.

In the case of the Dirichlet series for the Riemann zeta function, we
cannot use ansatz (\ref{Ansatz_RemSer}). This follows at once from the
relationship
\begin{equation}
%  \label{}
\Delta n^{\alpha} \; = \; (n+1)^{\alpha} - n^{\alpha} \; = \;
\alpha \, n^{\alpha-1} + \mathrm{O} \bigl( n^{\alpha-2} \bigr) \, ,
\qquad n \to \infty \, .
\end{equation}
Thus, we make the following ansatz, which takes into account the specific
features of the Dirichlet series (\ref{DiriSerZetaFun}):
\begin{equation}
  \label{ModRemSerZeta}
r_{n}^{(m)} \; - \, (n+2)^{1-s} \, \sum_{\mu=0}^{m} \,
\frac{\gamma_{\mu}}{(n+2)^{\mu}} \, , \qquad m \in \mathbb{N} \, , 
\quad n \in \mathbb{N}_0 \, .
\end{equation}
This ansatz is inspired by the truncation error estimate for Levin's $u$
transformation \cite{Levin/1973} (see also \cite[Section
7.3]{Weniger/1989} or \cite[Section IV]{Weniger/2004}).

As in Section \ref{Sec:AsyTruncErrAppr}, the unspecified coefficients
$\gamma_{\mu}$ are chosen in such a way that only a higher order
error remains:
\begin{equation}
  \label{ModRem_SerZeta_ConvCond}
\Delta r_{n}^{(m)} \; = \; (n+2)^{-s} \,
\bigl[ 1 \, + \, \mathrm{O} \bigl( n^{-m-1} \bigr) \bigr] \, ,
\qquad n \to \infty \, .
\end{equation}
For that purpose, we write:
\begin{align}
  \label{ModRem_SerZeta_ModExpr_1}
  \frac{r_{n+1}^{(m)} - r_{n}^{(m)}} {(n+2)^{1-s}} & = \;
  \sum_{\mu=0}^{m} \, \frac{\gamma_{\mu}}{(n+2)^{\mu}} \, - \, \left[
    \frac {n+3} {n+2} \right]^{1-s} \, \sum_{\mu=0}^{m} \,
  \frac{\gamma_{\mu}}{(n+3)^{\mu}}
  \\
  \label{ModRem_SerZeta_ModExpr_3}
  & = \; \sum_{\mu=0}^{m} \, \frac{\gamma_{\mu}}{(n+2)^{\mu}} \, \left\{
    1 \, - \, [1+1/(n+2)]^{1-s-\mu} \right\} \, .
\end{align}
With the help of the binomial series \cite[p.\
38]{Magnus/Oberhettinger/Soni/1966}, we obtain:
\begin{equation}
  \label{ModRem_SerZeta_InnerBraces}
1 - [1+1/(n+2)]^{1-s-\mu} \; = \; \sum_{\lambda=0}^{\infty} \,
\frac{(s+\mu-1)_{\lambda+1}}{(\lambda+1)!} \,
\frac{(-1)^\lambda}{(n+2)^{\lambda+1}} \, .
\end{equation}
Inserting (\ref{ModRem_SerZeta_InnerBraces}) into
(\ref{ModRem_SerZeta_ModExpr_3}) yields:
\begin{align}
%  \label{}
  & \sum_{\mu=0}^{m} \, \frac{\gamma_{\mu}}{(n+2)^{\mu}} \left\{ 1 \, -
    \, [1+1/(n+2)]^{1-s-\mu} \right\} \; = \;
  \notag \\
%  \label{}
  & \qquad \; = \; \sum_{\mu=0}^{m} \, \frac{\gamma_{\mu}}{(n+2)^{\mu}}
  \, \sum_{\lambda=0}^{\infty} \,
  \frac{(s+\mu-1)_{\lambda+1}}{(\lambda+1)!} \,
  \frac{(-1)^\lambda}{(n+2)^{\lambda+1}}
  \\
  & \qquad \; = \; - \, \sum_{\nu=0}^{\infty} \,(n+2)^{-\nu-1} \,
  \sum_{\lambda=0}^{\min (\nu, m)} \, \frac{(1-s-\nu)_{\lambda+1}
    \gamma_{\nu-\lambda}} { (\lambda+1)!} \, .
\end{align}
Thus, we obtain the following truncated asymptotic expansion:
\begin{align}
%  \label{}
  & \frac{r_{n+1}^{(m)} - r_{n}^{(m)}} {(n+2)^{-s}} \; = \; - \,
  \sum_{\mu=0}^{m} \,(n+2)^{-\mu}
  \notag \\
  & \qquad \times \, \sum_{\lambda=0}^{\mu} \,
  \frac{(1-s-\mu)_{\lambda+1} \gamma_{\mu-\lambda}}{(\lambda+1)!}  \, +
  \, \mathrm{O} \bigl(n^{-m-1}\bigr) \, , \qquad n \to \infty \, .
\end{align}
The unspecified coefficients $\gamma_{\mu}$ have to be determined
by solving the following system of linear equations, whose triangular
structure is obvious:
\begin{equation}
  \label{SysLinEqs_DirSerZeta} 
\sum_{\lambda=0}^{\mu} \, \frac{(1-s-\mu)_{\lambda+1}
\gamma_{\mu-\lambda}}{(\lambda+1)!} \; = \; \delta_{\mu 0} \, , 
\qquad 0 \le \mu \le m \, .
\end{equation}
For a more detailed analysis of the linear system
(\ref{SysLinEqs_DirSerZeta}), let us define $\beta_{\mu}$ via
\begin{equation}
  \label{Def_beta}
\gamma_{\mu} \; = \;
(-1)^{\mu} \, \frac{(s)_{\mu-1}}{{\mu}!} \, \beta_{\mu} \, ,
\qquad \mu \in \mathbb{N}_0 \, .
\end{equation}
Inserting (\ref{Def_beta}) into (\ref{SysLinEqs_DirSerZeta}) yields:
\begin{equation}
  \label{LinSys_beta_1}
\sum_{\lambda=0}^{\mu} \,
\frac{(1-s-\mu)_{\lambda+1}}{(\lambda+1)!} \,
\frac{(-1)^{\mu-\lambda} (s)_{\mu-\lambda-1}}{(\mu-\lambda)!} \,
\beta_{\mu-\lambda} \; = \; \delta_{\mu 0} \, , 
\qquad \mu \in \mathbb{N}_0 \, .
\end{equation}
Next, we use $(1-s-\mu)_{\lambda+1} = (-1)^{\lambda+1} (s)_{\mu}$ and
obtain 
\begin{align}
%  \label{}
  & \sum_{\lambda=0}^{\mu} \, \frac{(1-s-\mu)_{\lambda+1}
    \gamma_{\mu-\lambda}}{(\lambda+1)!}  \; = \; (-1)^{\mu+1} \,
  (s)_{\mu} \, \sum_{\lambda=0}^{\mu} \,
  \frac{\beta_{\mu-\lambda}}{(\lambda+1)! (\mu-\lambda)!}
  \\
%  \label{}
  & \qquad \; = \; \frac{(-1)^{\mu+1} (s)_{\mu}}{(\mu+1)!}  \,
  \sum_{\sigma=0}^{\mu} \, \frac{(\mu+1)!}{(\mu-\sigma+1)! {\sigma}!} \,
  \beta_{\sigma}
  \\
%  \label{}
  & \qquad \; = \; \frac{(-1)^{\mu+1} (s)_{\mu}}{(\mu+1)!}  \,
  \sum_{\sigma=0}^{\mu} \, \binom{\mu+1}{\sigma} \, \beta_{\sigma} \; =
  \; \delta_{\mu 0} \, , \qquad \mu \in \mathbb{N}_0 \, .
\end{align}
Thus, the linear system (\ref{SysLinEqs_DirSerZeta}) is equivalent to the
well known recurrence formula 
\begin{equation}
  \label{Rec_BN}
\sum_{\nu=0}^{n} \, \binom{n+1}{\nu} \, B_{\nu} \; = \; 0 \, ,
\qquad n \in \mathbb{N} \, ,
\end{equation} 
of the Bernoulli numbers (see for example \cite[Eq.\
(1.11)]{Temme/1996a}) together with the initial condition $B_0 = 1$.
Thus. the ansatz (\ref{ModRemSerZeta}) reproduces the finite
sum (\ref{EuMacZeta_a}) or (\ref{EuMacZeta_mod2_a}) of the
Euler-Maclaurin formula for the truncation error of the Dirichlet series,
which is not really surprising since asymptotic series are unique, if
they exist. Only the integral (\ref{EuMacZeta_b}) or
(\ref{EuMacZeta_mod2_b})) cannot be reproduced in this way.
 
\typeout{==> Section 6: The Gaussian Hypergeometric Series}
\section{The Gaussian Hypergeometric Series}
\label{Sec:Gauss_Hyg_2F1}

The simplicity of the terms of the Dirichlet series
(\ref{DiriSerZetaFun}) facilitates the derivation of explicit asymptotic
approximations to truncation errors by solving a system of linear
equations in closed form. A much more demanding test for the feasibility
of the new formalism is the Gaussian hypergeometric series
(\ref{Ser_2F1}), which depends on three parameters $a$, $b$, and $c$, and
one argument $z$.

Due to the complexity of the terms of the Gaussian hypergeometric series
(\ref{Ser_2F1}), there is little hope in obtaining explicit analytical
solutions to the linear equations. From a pragmatist's point of view, it
is therefore recommendable to use computer algebra systems like Maple and
Mathematica and let the computer do the work.

In the case of a nonterminating Gaussian hypergeometric series, we have:
\begin{align}
 \label{ParSum_2F1}
 s_n (z) & \; = \; \sum_{\nu=0}^{n} \, \frac {(a)_{\nu} (b)_{\nu}}
 {(c)_{\nu} {\nu}!} \, z^{\nu} \, ,
 \\
 \label{Rem_2F1}
 r_n (z) & \; = \; - \sum_{\nu=n+1}^{\infty} \, \frac {(a)_{\nu}
   (b)_{\nu}} {(c)_{\nu} {\nu}!} \, z^{\nu}
 \notag \\
 & \; = \; - \frac{(a)_{n+1} (b)_{n+1}}{(c)_{n+1} (n+1)!} z^{n+1} \,
 \sum_{\nu=0}^{\infty} \, \frac {(a+n+1)_{\nu} (b+n+1)_{\nu}}
 {(c+n+1)_{\nu} (n+2)_{\nu}} \, z^{\nu} \, ,
 \\
  \label{FinDifEq_2F1}
  \Delta r_n (z) & \; = \; \frac {(a)_{n+1} (b)_{n+1}} {(c)_{n+1}
    (n+1)!} \, z^{n+1} \, .
\end{align}

Since $[(a+n+1)_{\nu} (b+n+1)_{\nu}]/[(c+n+1)_{\nu} (n+2)_{\nu}]$ can be
expressed as a power series in $1/(n+1)$, the following ansatz make
sense:
\begin{align}
  \label{ModRem_2F1}
  r_{n}^{(m)} (z) & \; = \; - \, \frac {(a)_{n+1} (b)_{n+1}} {(c)_{n+1}
    (n+1)!} \, z^{n+1} \, \sum_{\mu=0}^{m} \,
  \frac{\gamma_{\mu}}{(n+1)^{\mu}} \, ,
  \notag \\
  & \qquad m \in \mathbb{N} \, , \quad n \in \mathbb{N}_0 \, , \quad
  \vert z \vert < 1 \, .
\end{align}
Again, we choose the unspecified coefficients $\gamma_{\mu}$ in
(\ref{ModRem_2F1}) in such a way that only a higher order error error
remains:
\begin{align}
  \label{ModRem_2F1_ConvCond}
  & \Delta r_{n}^{(m)} (z) \; = \; r_{n+1}^{(m)} (z) - r_{n}^{(m)} (z)
  \notag \\
  & \qquad \; = \; \frac {(a)_{n+1} (b)_{n+1}} {(c)_{n+1} (n+1)!} \,
  z^{n+1} \, \bigl[ 1 \, + \, \mathrm{O} \bigl( n^{-m-1} \bigr) \bigr]
  \, , \qquad n \to \infty \, .
\end{align}
This convergence condition can be reformulated as follows:
\begin{align}
  & \frac{r_{n+1}^{(m)} (z) - r_{n}^{(m)} (z)} {[(a)_{n+1} (b)_{n+1}
    z^{n+1}]/[(c)_{n+1} (n+1)!]}
  \notag \\
  \label{ModRem_2F1_ModExpr}
  & \qquad = \; \sum_{\mu=0}^{m} \, \frac{\gamma_{\mu}}{(n+1)^{\mu}} \,
  - \, \frac {(a+n+1) (b+n+1)} {(c+n+1) (n+2)} \, z \, \sum_{\mu=0}^{m}
  \, \frac{\gamma_{\mu}}{(n+2)^{\mu}}
  \\
  \label{ModRem_2F1_ModConvCond}
  & \qquad = \; 1 \, + \, \mathrm{O} \bigl( n^{-m-1} \bigr) \, , \qquad
  n \to \infty \, .
\end{align}
Now, we only have to do an asymptotic expansion of
(\ref{ModRem_2F1_ModExpr}) in terms of the asymptotic sequence $\{
1/(n+1)^{j} \}_{j=0}^{\infty}$ as $n \to \infty$. This yields:
\begin{align}
  \label{AsyRat_2F1}
  & \frac{r_{n+1}^{(m)} (z) - r_{n}^{(m)} (z)} {[(a)_{n+1} (b)_{n+1}
    z^{n+1}]/[(c)_{n+1} (n+1)!]}
  \notag \\
%  \label{}
  & \qquad = \; \sum_{\mu=0}^{m} \,
  \frac{\mathcal{C}_{\mu}}{(n+1)^{\mu}} \, + \, \mathrm{O} \bigl(
  n^{-m-1} \bigr) \, , \qquad n \to \infty \, .
\end{align}
We then obtain the following system of coupled linear equations in the
unspecified coefficient $\gamma_{\mu}$ with $0 \le \mu \le m$:
\begin{equation}
  \label{LinSys_2F1}
\mathcal{C}_{\mu} \; = \; \delta_{\mu 0} \, , 
\qquad 0 \le \mu \le m \, .
\end{equation}

As discussed in Section \ref{Sec:AsyTruncErrAppr}, a coefficient
$\mathcal{C}_{\mu}$ with $0 \le \mu \le m$ contains only the unspecified
coefficients $\gamma_{0}$, $\dots$, $\gamma_{\mu}$ but not
$\gamma_{\mu+1}$ $\dots$, $\gamma_{m}$. Thus, the symbolic solution of
these linear equations for a Gaussian hypergeometric function ${}_2 F_1
(a, b; c; z)$ with unspecified parameters $a$, $b$, and $c$ and
unspecified argument $z$ is not particularly difficult for a computer
algebra system, since the unspecified coefficient $\gamma_{\mu}$ can be
determined successively. The following linear equations were constructed
with the help of Maple 8:
\begin{subequations}
   \label{LinSys_2F1_expl}
  \begin{align}
%  \label{}
    \mathcal{C}_0 & \; = \; (1 - z) \, \gamma_{0} \; = \; 1 \, ,
    \\
%  \label{}
    \mathcal{C}_1 & \; = \; ( c - a - b + 1) \, z \, \gamma_{0} + (1 -
    z) \, \gamma_{1} \; = \; 0 \, ,
    \\
%  \label{}
    \mathcal{C}_2 & \; = \; [(c - b + 1) \, a + (c + 1) \, b - 1 - c -
    c^{2}] \, z \, \gamma_{0}
    \notag \\
    & \qquad + \, (c + 2 - b - a) \, z \, \gamma_{1} + (1 - z) \,
    \gamma_{2} \; = \; 0 \, ,
    \\
%  \label{}
    \mathcal{C}_3 & \; = \; \{[(c + 1)\,b - 1 - c - c^{2}] \, a - (1 + c
    + c^{2}) \, b + c^{3} + c^{2} + c + 1 \} \, z \, \gamma_{0}
    \notag \\
    & \qquad + \, [(c + 2 - b)\,a + (c + 2)\,b - 3 - c^{2} - 2\,c] \, z
    \, \gamma_{1}
    \notag \\
    & \qquad + \, (3 - b - a + c) \, z \, \gamma_{2} + (1 - z) \,
    \gamma_{3} \; = \; 0 \, .
  \end{align}  
\end{subequations} 
This example shows that the complexity of the coefficients
$\mathcal{C}_{\mu}$ in (\ref{AsyRat_2F1}) increases so rapidly with
increasing index $\mu$ that a solution of the linear equations
(\ref{LinSys_2F1}) becomes soon unmanageable for humans. This is also
confirmed by the following solutions of (\ref{LinSys_2F1_expl}) obtained
symbolically with the help of Maple 8:
\begin{subequations}
   \label{SolSys_2F1_expl}
\begin{align}
%  \label{}
(z-1) \, \gamma_{0} & \; = \; 1 \, ,
\\
%  \label{}
(z-1)^2 \, \gamma_{1} & \; = \; z \, (a + b - c - 1) \, ,
\\
%  \label{}
(z-1)^3 \, \gamma_{2} & \; = \; z \, \{ [a^{2} + (b - c - 2) \, a
+ b^{2} - (c + 2) \, b + 1 + 2 \, c] \, z
\notag \\
& \qquad + (b - c - 1) \, a - ( c + 1)\,b + 1 + c + c^{2}\} \, ,
\\
%  \label{}
(z-1)^4 \, \gamma_{3} & \; = \; z \, \{ [a^{3} +
(b - c - 3) \, a^{2} + (b^{2} - (c + 3) \,b + 3 + 3 \, c) \, a
\notag \\
& \qquad
+ b^{3} - (c + 3)\,b^{2} + (3 + 3 \, c) \, b - 1 - 3 \, c] \, z^{2}
\notag \\
& \qquad
+ [(2 \, b - 2 \, c - 3) \, a^{2} + (2 \, b^{2} - (4 \, c + 8) \, b
+ 2 \, c^{2} + 7 + 8\,c)\,a
\notag \\
& \qquad - (2 \, c + 3) \, b^{2} + (2\,c^{2}
 + 7 + 8\,c)\,b - 4 - 5\,c^{2} - 7 \, c] \, z
\notag \\
& \qquad + [- (c + 1) \, b + 1 + c^{2} + c] \, a + (1 + c^{2} + c) \, b
\notag \\
& \qquad - 1 - c^{2} - c - c^{3} \} \, .
\end{align}  
\end{subequations}

The solutions (\ref{SolSys_2F1_expl}), which are rational in $z$,
demonstrate quite clearly a principal weakness of symbolic computing.
Typically, the results are complicated and poorly structured algebraic
expressions, and it is normally very difficult to gain further insight
from them. Nevertheless, symbolic solutions of the linear equations
(\ref{LinSys_2F1}) are computationally very useful.

For the Gaussian hypergeometric series with $a = 1/3$, $b = 7/5$, $c =
9/2$, and $z = - 0.85$, Maple~8 produced for $m = 8$ and $n = 1$ the
following results:
\begin{subequations}
%  \label{}
\begin{align}
%  \label{}
r_{1} & \; = \; - 0.016~412~471 \, ,
\\
%  \label{}
a_{2} \, [4/4] & \; = \; - 0.016~410~482 \, ,
\\
%  \label{}
a_{2} \, \tilde{r}_{1}^{(8)} & \; = \; - 0.016~414~203 \, ,
\\
%  \label{}
a_{2} \, r_{1}^{(8)} & \; = \; - 0.004~008~195 \, .
\end{align}
\end{subequations}
It is in my opinion quite remarkable that for $n = 1$, which is very far
away from the asymptotic regime, at least the Pad\'{e} approximant
$a_{2} [4/4]$ and the truncated factorial series $a_{2}
\tilde{r}_{1}^{(8)}$ agree remarkably well with the ``exact'' truncation
error $r_1$. In contrast, the truncated inverse power series $a_{2}
r_{1}^{(8)}$ produces a relatively poor result. For $n = 10$, which
possibly already belongs to the asymptotic regime, Maple~8 produced the
following results:
\begin{subequations}
%  \label{}
\begin{align}
%  \label{}
r_{10} & \; = \; 0.000~031~925~482 \, ,
\\
%  \label{}
a_{11} \, [4/4] & \; = \; 0.000~031~925~482 \, ,
\\
%  \label{}
a_{11} \, \tilde{r}_{10}^{(8)} & \; = \; 0.000~031~925~483 \, ,
\\
%  \label{}
a_{11} \, r_{10}^{(8)} & \; = \; 0.000~031~925~471 \, .
\end{align}
\end{subequations}

Finally, let me emphasize that the formalism of this article is not
limited to a Gaussian hypergeometric series (\ref{Ser_2F1}), but works
just as well in the case of a generalized hypergeometric series ${}_{p+1}
F_p \bigl( \alpha_1, \dots, \alpha_{p+1}.  \beta_1, \dots, \beta_p; z
\bigr)$.

\typeout{==> Section 7: The Asymptotic Series for the Exponential Integral}
\section{The Asymptotic Series for the Exponential Integral}
\label{Sec:AsySer_E_1}

The divergent asymptotic series (\ref{AsySerE_1}) for the exponential
integral $E_1 (z)$ is probably the most simple model for many other
factorially divergent asymptotic inverse power series occurring in
special function theory.  Well known examples are the asymptotic series
for the modified Bessel function $K_{\nu} (z)$, the complementary error
function $\mathrm{erfc} (z)$, the incomplete gamma function $\Gamma (a,
z)$, or the Whittaker function $W_{\kappa,\mu} (z)$.  Moreover,
factorial divergence is also the rule rather than the exception among
the perturbation expansions of quantum physics (see \cite{Weniger/2004}
for a condensed review of the relevant literature).

The exponential integral $E_1 (z)$ can also be expressed as a Stieltjes
integral:
\begin{equation}
  \label{StieInt_E_1}
z \, \mathrm{e}^z \, E_1 (z) \; = \; \int_{0}^{\infty} \, 
\frac{\mathrm{e}^{-t} \mathrm{d} t}{1+t/z} \, .
\end{equation}
If $z < 0$, this integral has to be interpreted as a principal value
integral.

In the case of a factorially divergent inverse power series, it is of
little use to represent the truncation error $r_n (z)$ by a power series
as in (\ref{Rem_2F1}). If, however, we use $\sum_{\nu=0}^{n} x^{\nu} =
[1-x^{n+1}]/[1-x]$ in (\ref{StieInt_E_1}), we immediately obtain:
\begin{align}
 \label{ParSum_E_1}
 s_n (z) & \; = \; \sum_{\nu=0}^{n} \, (-1/z)^{\nu} \, \nu! \, ,
 \\
 \label{Rem_E_1}
 r_n (z) & \; = \; - (-z)^{-n-1}
\, \int_{0}^{\infty} \,
\frac{t^{n+1} \mathrm{e}^{-t} \mathrm{d}t}{1 + t/z} \, ,
 \\[1.5\jot]
  \label{FinDifEq_E_1}
  \Delta r_n (z) & \; = \; (-1/z)^{n+1} \, (n+1)! \, .
\end{align}
Because of the factorial growth of the coefficients in
(\ref{AsySerE_1}), it is advantageous to use instead of
(\ref{Ansatz_RemSer}) the following ansatz:
\begin{equation}
  \label{ModRem_E_1}
r_{n}^{(m)} (z) \; = \; (-1/z)^n n! \, \sum_{\mu=0}^{m} \,
  \frac{\gamma_{\mu}}{(n+1)^{\mu}} \, ,
\qquad m \in \mathbb{N} \, , \quad n \in \mathbb{N}_0 \, .
\end{equation}
Again, we choose the unspecified coefficients $\gamma_{\mu}$ in
(\ref{ModRem_E_1}) in such a way that only a higher order error error
remains:
\begin{equation}
  \label{ModRem_E_1_ConvCond}
\Delta r_{n}^{(m)} (z) \; = \; (-1/z)^{n+1} \, (n+1)! \, 
   \bigl[ 1 \, + \, \mathrm{O} \bigl( n^{-m-1} \bigr) \bigr] \,
   , \qquad n \to \infty \, .
\end{equation}
This convergence condition can be reformulated as follows:
\begin{align}
  \label{ModRem_E_1_ModExpr}
  \frac{r_{n+1}^{(m)} (z) - r_{n}^{(m)} (z)} {(-1/z)^{n+1} (n+1)!}
  & \; = \; \frac{-z}{n+1} \, \sum_{\mu=0}^{m} \,
  \frac{\gamma_{\mu}}{(n+1)^{\mu}} \, - \,  \sum_{\mu=0}^{m} \,
  \frac{\gamma_{\mu}}{(n+2)^{\mu}}
  \\
  \label{ModRem_E_1_ModConvCond}
  & \; = \; 1 \, + \, \mathrm{O} \bigl( n^{-m-1} \bigr) \, ,
  \qquad n \to \infty \, .
\end{align}
Next, we do an asymptotic expansion of the right-hand side of
(\ref{ModRem_E_1_ModExpr}) in terms of the asymptotic sequence $\{
1/(n+1)^{j} \}_{j=0}^{\infty}$ as $n \to \infty$. This yields:
\begin{equation}
  \label{AsyRat_E_1}
\frac{r_{n+1}^{(m)} (z) - r_{n}^{(m)} (z)}
{(-1/z)^{n+1} (n+1)!} \; = \;
\sum_{\mu=0}^{m} \, \frac{\mathcal{C}_{\mu}}{(n+1)^{\mu}} \, + \, 
\mathrm{O} \bigl( n^{-m-1} \bigr) \, , \qquad n \to \infty \, .
\end{equation}
Again, we have to solve the following system of linear equations:
\begin{equation}
  \label{LinSys_E_1}
\mathcal{C}_{\mu} \; = \; \delta_{\mu 0} \, , 
\qquad 0 \le \mu \le m \, .
\end{equation}
The following linear equations were constructed with the help of Maple 8:
\begin{subequations}
  \label{LinSys_E_1_expl}
\begin{align}
%  \label{}
\mathcal{C}_0 & \; = \; - \gamma_{0} \; = \; 1 \, ,
\\
%  \label{}
\mathcal{C}_1 & \; = \; - z \, \gamma_{0} - \gamma_{1} \; = \; 0 \, ,
\\
%  \label{}
\mathcal{C}_2 & \; = \; (1 - z) \, \gamma_{1} - \gamma_{2} \; = \; 0 \, ,
\\
%  \label{}
\mathcal{C}_3 & \; = \; - \gamma_{1} + (2 - z) \, \gamma_{2} 
- \gamma_{3} \; = \; 0 \, ,
\\
%  \label{}
\mathcal{C}_4 & \; = \; \gamma_{1} - 3 \, \gamma_{2} + 
(3 - z) \, \gamma_{3} - \gamma_{4} \; = \; 0 \, .
\end{align}
\end{subequations}
If we compare the complexity of the equations (\ref{LinSys_E_1_expl})
with those of (\ref{LinSys_2F1_expl}), we see that in the case of the
asymptotic series (\ref{AsySerE_1}) for the exponential integral there
may be a chance of finding explicit expressions for the coefficients
$\gamma_{\mu}$. At least, the solutions of the linear system
(\ref{LinSys_E_1_expl}) obtained symbolically with the help of Maple 8
look comparatively simple:
\begin{subequations}
  \label{SolSys_E_1_expl}
\begin{align}
%  \label{}
\gamma_{0} & \; = \; - 1 \, ,
\\
%  \label{}
\gamma_{1} & \; = \; z \, ,
\\
%  \label{}
\gamma_{2} & \; = \; - (z - 1) \, z \, ,
\\
%  \label{}
\gamma_{3} & \; = \; (z^{2} - 3 \, z + 1) \, z \, ,
\\
%  \label{}
\gamma_{4} & \; = \; - (z^{3} - 6 \, z^{2} + 7 \, z - 1) \,z \, .
\end{align}
\end{subequations}
Of course, this requires further investigations.

The relative simplicity of the coefficients in (\ref{SolSys_E_1_expl})
offers other perspectives. For example, the $[2/2]$ Pad\'{e} approximant
to the truncated power series in (\ref{ModRem_E_1}) is compact enough to
be printed without problems:
\begin{equation}
%  \label{Pade[2/2]_ExpInt}
[2/2] \; = \; \frac
{\displaystyle - 1 + \frac{3 - z}{n + 1} + \frac{2}{(n + 1)^{2}}}
{\displaystyle 1 + \frac{2z - 3}{n + 1} + 
\frac{z^{2} - 2z + 2}{(n + 1)^{2}}}
\; = \;
\frac{n^{2} - n + zn + z}{n^{2} - n + 2zn + z^{2}} \, .
\end{equation}
Pad\'{e} approximants to the truncated power series in (\ref{ModRem_E_1})
seem to be a new class of approximants that are rational in both $n$ and
$z$.

For the asymptotic series (\ref{AsySerE_1}) for $E_1 (z)$ with $z = 5$,
Maple~8 produced for $m = 16$ and $n = 2$ the following results:
\begin{subequations}
%  \label{}
\begin{align}
%  \label{}
r_{2} & \; = \; 0.027~889 \, ,
\\
%  \label{}
a_{2} \, [8/8] & \; = \; 0.027~965 \, ,
\\
%  \label{}
a_{2} \, \tilde{r}_{2}^{(16)} & \; = \; 0.028~358 \, ,
\\
%  \label{}
a_{2} \, r_{2}^{(16)} & \; = \; - 177.788 \, .
\end{align}
\end{subequations}
The Pad\'{e} approximant $a_{2} [8/8]$ and the truncated factorial series
$a_{2} \tilde{r}_{2}^{(16)}$ agree well with the ``exact'' truncation
error $r_2$, but the truncated inverse power series $a_{2} r_{2}^{(16)}$
is way off. For $n = 10$, all results agree reasonably well
\begin{subequations}
%  \label{}
\begin{align}
%  \label{}
r_{10} & \; = \; 0.250~470~879 \, ,
\\
%  \label{}
a_{10} \, [8/8] & \; = \; 0.250~470~882 \, ,
\\
%  \label{}
a_{10} \, \tilde{r}_{10}^{(16)} & \; = \; 0.250~470~902 \, ,
\\
%  \label{}
a_{10} \, r_{10}^{(16)} & \; = \; 0.250~470~221 \, .
\end{align}
\end{subequations}

\typeout{==> Section 8: Summary and Outlook}
\section{Conclusions and Outlook}
\label{Sec:ConclusionsOutlook}

A new formalism is proposed that permits the construction of asymptotic
approximations to truncation errors $r_n = - \sum_{\nu=n+1}^{\infty}
a_{\nu}$ of infinite series for special functions by solving a system of
linear equations. Approximations to truncation errors of monotone series
can be obtained via the Euler-Maclaurin formula. The formalism proposed
here is, however, based on different assumptions and can be applied even
if the terms of the series have a comparatively complicated structure. In
addition, the new formalism works also in the case of alternating and
even divergent series.

Structurally, the asymptotic approximations of this article resemble the
asymptotic inverse power series for special functions as $z \to \infty$,
since they are not expansions of $r_n$, but rather expansions of ratios
like $r_n/a_{n+1}$, $r_n/a_{n}$, or $r_n/[(n+\alpha) a_{n+1}]$ with
$\alpha > 0$. This is consequential, because it makes it possible to use
the convenient asymptotic sequence $\{ 1/(n+\alpha)^{j}
\}_{j=0}^{\infty}$ of inverse powers. This greatly facilitate the
necessary analytical manipulations and ultimately leads to comparatively
simple systems of linear equations.

As shown in Section \ref{Sec:DirSerRiemannZetaFun}, the new formalism
reproduces in the case of the Dirichlet series (\ref{DiriSerZetaFun}) for
the Riemann zeta function the expressions (\ref{EuMacZeta_a}) or
(\ref{EuMacZeta_mod2_a}) that follow from the Euler-Maclaurin formula.
The linear equations (\ref{SysLinEqs_DirSerZeta}) are equivalent to the
recurrence formula (\ref{Rec_BN}) of the Bernoulli numbers. Thus, only
the integral (\ref{EuMacZeta_b}) or (\ref{EuMacZeta_mod2_b})) cannot be
obtained in this way.

Much more demanding is the Gaussian hypergeometric series
(\ref{Ser_2F1}), which is discussed in Section \ref{Sec:Gauss_Hyg_2F1}.
The terms of this series depend on three in general complex parameters
$a$, $b$, and $c$ and one argument $z$. Accordingly, there is little hope
that we might succeed in finding an explicit solution to the linear
equations.  However, all linear equations considered in this article have
a triangular structure. Consequently, it is relatively easy to construct
solutions symbolically with the help of a computer algebra system like
Maple. The numerical results presented in Section \ref{Sec:Gauss_Hyg_2F1}
also indicate that the formalism proposed in this article is indeed
computationally useful.

As a further example, the divergent asymptotic series (\ref{AsySerE_1})
for the exponential integral $E_1 (z)$ is considered in Section
\ref{Sec:AsySer_E_1}. The linear equations are again solved symbolically
by Maple. Numerical results are also presented. This example is important
since it shows that the new formalism works also in the case of
factorially divergent series. The Euler-Maclaurin formula can only handle
convergent monotone series.

Although the preliminary results look encouraging, a definite assessment
of the usefulness of the new formalism for the computation of special
functions is not yet possible. This requires much more data.
Consequently, the new formalism should be be applied to other series
expansions for special functions and the performance of the resulting
approximations should be analyzed and compared with other computational
approaches.

I suspect that in most cases it will be necessary to solve the linear
equations symbolically with the help of a computer algebra system like
Maple. Nevertheless, it cannot be ruled out that at least for some
special functions with sufficiently simple series expansions explicit
analytical solutions to the linear equations can be found.

Effective numerical analytic continuation methods are of considerable
relevance for the new formalism which produces asymptotic approximations.
We cannot tacitly assume that these approximations provide good results
outside the asymptotic regime, although it would be highly desirable to
use them also for small indices. In Section \ref{Sec:NumAnaCon}, only
factorial series and Pad\'{e} approximants are mentioned, although many
other numerical techniques are known that can accomplish such an analytic
continuation. Good candidates are sequence transformations which are
often more effective than the better known Pad\'{e} approximants.
Details can be found in books by Brezinski and Redivo Zaglia
\cite{Brezinski/RedivoZaglia/1991a}, Sidi \cite{Sidi/2003}, or Wimp
\cite{Wimp/1981}, or in a review by the present author
\cite{Weniger/1989}.

%%%%%%%%%%%%%%%%%%%%%%%%%%%%%%%%%%%%%%%%%%%%%%%%%%%%%%%%%%%%%%%%%%%%%%
%
% References
%

\providecommand{\SortNoop}[1]{} \providecommand{\OneLetter}[1]{#1}
  \providecommand{\SwapArgs}[2]{#2#1}

%
%%%%%%%%%%%%%%%%%%%%%%%%%%%%%%%%%%%%%%%%%%%%%%%%%%%%%%%%%%%%%%%%%%%%%%

%\printindex
\end{document}